\documentclass[12pt]{amsart} 

\usepackage{amsmath} \usepackage{amssymb} \usepackage{amsthm}
\usepackage{amscd} 

\def\a{{\mathfrak{a}}} \def\b{{\mathfrak{b}}}
\def\F{{\mathbb{F}}} \def\J{{\mathcal{J}}} \def\m{{\mathfrak{m}}} \def\Z{{\mathbb{Z}}}\def\N{{\mathbb{N}}} 
\def\Q{{\mathbb{Q}}}  
\def\Hom{{\mathrm{Hom}}}  
  \def\Ann{{\mathrm{Ann}}}
\def\Spec{{\mathrm{Spec\; }}}

\theoremstyle{plain}
\newtheorem{thm}{Theorem}[section] 
\newtheorem*{mainthm}{Main Theorem}
\newtheorem{cor}[thm]{Corollary}
\newtheorem{prop}[thm]{Proposition}
\newtheorem{conj}[thm]{Conjecture}
\newtheorem{propdef}[thm]{Proposition-Definition} 

\newtheorem{lem}[thm]{Lemma}
\theoremstyle{definition} 
\newtheorem{defn}[thm]{Definition}
 
\theoremstyle{remark}
\newtheorem{rem}[thm]{Remark}

\newtheorem*{cl}{Claim}
\newtheorem{cln}{Claim}
\newtheorem*{acknowledgement}{Acknowledgments}

\title{Generalized test ideals and symbolic powers}
\author{Shunsuke Takagi}
\address{Department of Mathematics, Kyushu University, 
6-10-1 Hakozaki, Higashi-ku, Fukuoka, 812-8581 JAPAN}
\email{stakagi@math.kyushu-u.ac.jp}
\author{Ken-ichi Yoshida}
\address{Graduate School of Mathematics, Nagoya University, Chikusa-ku, Nagoya,464-8602 JAPAN}
\email{yoshida@math.nagoya-u.ac.jp}

\begin{document}
\maketitle
\markboth{S.TAKAGI and K.YOSHIDA}{Generalized test ideals and Symbolic powers}

\begin{abstract}
Hochster and Huneke proved in \cite{HH5}  fine behaviors of symbolic powers of ideals in regular rings, using the theory of tight closure. In this paper, we use generalized test ideals, which are a characteristic $p$ analogue of multiplier ideals, to give a slight generalization of Hochster-Huneke's results. 
\end{abstract}

\section*{Introduction}
Ein, Lazarsfeld and Smith proved in \cite{ELS} the following uniform behavior of symbolic powers of ideals in affine regular rings of equal characteristic zero: if $h$ is the largest height of any associate prime of an ideal $I \subseteq R$, then $I^{(hn+kn)} \subseteq (I^{(k+1)})^n$ for all integers $n \ge 1$ and $k \ge 0$. Here, if $W$ is the complement of the union of the associate primes of $I$, then $m^{\rm th}$ symbolic powers $I^{(m)}$ of $I$ is defined to be the contraction of $I^mR_W$ to $R$, where $R_W$ is the localization of $R$ at the multiplicative system $W$. 
To prove this, they introduced the notion of asymptotic multiplier ideals, which is a variant of multiplier ideals associated to filtrations of ideals and is formulated in terms of resolution of singularities. 
The above uniform behavior of symbolic powers immediately follows from a combination of properties of asymptotic multiplier ideals whose proofs need deep vanishing theorems. 

After a short while, Hochster and Huneke generalized in \cite{HH4} Ein-Lazarsfeld-Smith's result to the case of arbitrary regular rings of equal characteristic (i.e.,~both of equal characteristic zero and of positive prime  characteristic) in a completely different way. 
Furthermore, they used in \cite{HH5} similar ideas to prove more subtle behaviors of symbolic powers of ideals in a regular ring of equal characteristic.
Their methods depend on the theory of tight closure and reduction to positive characteristic, and, in consequence, they need neither resolution of singularities nor vanishing theorems which are proved only in characteristic zero. 
In this paper, combining the ideas of Ein-Lazarsfeld-Smith and Hochster-Huneke, we give a slight generalization of Hochster-Huneke's results in \cite{HH5}. 

Tight closure is an operation defined on ideals or modules in positive characteristic, introduced by Hochster-Huneke \cite{HH1} in the 1980s.
The test ideal $\tau(R)$ of a Noetherian ring $R$ of prime characteristic $p$ is the annihilator ideal of all tight closure relations in $R$, and it plays a central role in the theory of tight closure. In \cite{HY} and \cite{Ha}, Hara and Yoshida introduced a generalization of the test ideal $\tau(R)$, the ideal $\tau(\a_{\bullet})$ associated to a filtration of ideals $\a_{\bullet}$, and showed that their generalized test ideal $\tau(\a_{\bullet})$ is a characteristic $p$ analogue of the asymptotic multiplier ideal $\J(\a_{\bullet})$. In particular, in fixed prime characteristic, the generalized test ideal $\tau(\a_{\bullet})$ satisfies several nice properties similar to those of the asymptotic multiplier ideal $\J(\a_{\bullet})$ which are needed to prove Ein-Lazarsfeld-Smith's result; e.g., an analogue of Skoda's theorem (\cite[Theorems 4.1, 4.2]{HT}) and the subadditivity theorem (\cite[Theorem 4.5]{HY}, \cite[Proposition 2.10]{Ha}). 
In this paper, we first prove that the formation of generalized test ideals commutes with localization without the assumption of F-finiteness (compare this with \cite[Proposition 3.1]{HT}). 
This result makes easier to study further properties of generalized test ideals in non-F-finite rings.
Then, employing the strategy of  Ein-Lazarsfeld-Smith, we use generalized test ideals instead of asymptotic multiplier ideals to prove the following behavior of symbolic powers.

\begin{mainthm}
Let $R$ be an excellent regular ring of characteristic $p>0$ $($resp. regular algebra essentially of finite type over a field of characteristic zero$)$ and let $I \subsetneq R$ be an ideal of positive height. 
Let $h$ denote the largest analytic spread of $IR_P$ as $P$ runs through the associated primes of $I$.

\begin{enumerate}
\item 
If $(R,\m)$ is local, then for all integers $n \ge 1$ and $k \ge 0$, one has 
$$I^{(hn+kn+1)}\subseteq \m(I^{(k+1)})^n.$$
\item
If $R/I$ is F-pure $($resp. of dense F-pure type$)$,  then for all integers $n \ge 1$ and $k \ge 0$, one has
 $$I^{(hn+kn-1)} \subseteq (I^{(k+1)})^n.$$ 
\end{enumerate}
\end{mainthm}

These results are a slight generalization of the results in \cite{HH5} and are closely related to Eisenbud-Mazur's conjecture concerning the existence of evolutions:
Eisenbud and Mazur asked in \cite{EM} whether $P^{(2)} \subseteq \m P$ for a prime ideal $P$ in a regular local ring $(R,\m)$ of equal characteristic zero. 
Their conjecture fails in positive characteristic (see \cite{EM}), whereas our results hold true even in positive characteristic.  The authors hope that their methods shed a new light on the theory of symbolic powers of ideals. 

\begin{acknowledgement}
The first author thanks Craig Huneke for sending him his preprint \cite{HH5} prior to publication and for answering several questions. 
He is also indebted to Lawrence Ein and Sean Sather-Wagstaff for valuable conversation. 
The authors were partially supported by Grant-in-Aid for Scientific Research, 17740021 and 19340005, respectively, from JSPS. 
The first author was also partially supported by Program for Improvement of Research Environment for Young Researchers from SCF commissioned by MEXT of Japan.
\end{acknowledgement}

\section{Generalized test ideals}
In \cite{Ha}, Hara generalized the notion of tight closure to define a positive characteristic analogue of asymptotic multiplier ideals, which are ``multiplier ideals associated to graded families of ideals" (see below for the definition of graded families of ideals). 
In this section, we quickly review the definition and basic properties of this analogue. 

In this paper, all rings are excellent Noetherian reduced commutative rings with unity. 
For a ring $R$, we denote by $R^{\circ}$ the set of elements of $R$ which are not in any minimal prime ideal.  
A \textit{graded family of ideals} $\a_{\bullet}=\{ \a_m \}_{m \ge 1}$ on $R$ means a collection of ideals $\a_m \subseteq R$, satisfying $\a_1 \cap R^{\circ} \ne \emptyset$ and $\a_k \cdot \a_l \subseteq \a_{k+l}$ for all $k, l \ge 1$. 
Just for convenience, we decree that $\a_0=R$. 
One of the most important examples of graded families of ideals is a collection of symbolic powers $\a^{(\bullet)}=\{\a^{(m)}\}_{m \ge 1}$. 

Let $R$ be a ring of characteristic $p >0$ and $F\colon R \to R$ be the Frobenius map which sends $x \in R$ to $x^p \in R$. 
For each integer $e > 0$, the ring $R$ viewed as an $R$-module via the $e$-times iterated Frobenius map $F^e \colon R \to R$ is denoted by ${}^e\! R$. 
Since $R$ is assumed to be reduced, we can identify $F^e \colon R \to {}^e\! R$ with the natural inclusion map $R \hookrightarrow R^{1/p^e}$. 
We say that $R$ is \textit{F-pure} if $R^{1/p}$ is a pure extension of $R$ and that $R$ is {\it F-finite} if $R^{1/p}$ is a finitely generated $R$-module. For example, complete local rings with perfect residue fields are F-finite. 
An F-finite ring $R$ is said to be \textit{strongly F-regular} if for any $c \in R^{\circ}$, there exists $q=p^e$ such that $c^{1/q}R \hookrightarrow R^{1/q}$ splits as an $R$-linear map. 

Let $R$ be a ring of characteristic $p > 0$ and $M$ be an $R$-module. 
For each integer $e > 0$, we denote $\F^e(M) = \F_R^e(M) := M \otimes_R {}^e\! R$ and regard it as an $R$-module by the action of $R$ on ${}^e\! R$ from the right. 
Then we have the induced $e$-times iterated Frobenius map $F_M^e \colon M \to \F^e(M)$. 
The image of $z \in M$ via this map is denoted by $z^q:= F_M^e(z) 
\in \F^e(M)$, where $q=p^e$. For an $R$-submodule $N$ of $M$, we denote by $N^{[q]}_M$ the 
image of the induced map $\F^e(N) \to \F^e(M)$. 

Now we recall the definition of $\a_{\bullet}^k$-tight closure, which is a variant of tight closure associated  to a graded family of ideals $\a_{\bullet}$ with exponent $k$. 
\begin{defn}[\textup{\cite[Definition 2.7]{Ha}}]
Let $\a_{\bullet}=\{\a_m \}$ be a graded family of ideals on a ring $R$ of characteristic $p>0$ and let $k \ge 1$ be an integer. 
Let $N \subseteq M$ be (not necessarily finitely generated) $R$-modules. 
The $\a_{\bullet}^k$-tight closure of $N$ in $M$, denoted by $N_M^{*\a_{\bullet}^k}$, is defined to be the submodule of $M$ consisting of all elements $z \in M$ for which there exists $c \in R^{\circ}$ such that 
$$c \a_{kq} z^q \subseteq N_M^{[q]}$$
for all large $q=p^e$. 
The $\a_{\bullet}^{k}$-tight closure $I^{*\a_{\bullet}^k}$ of an ideal $I \subseteq R$ is just defined by $I^{*\a_{\bullet}^{k}} = I^{*\a_{\bullet}^{k} }_R$.
\end{defn} 

\begin{rem}
When $\a_m=R$ for all $m \ge 1$, ${\a_{\bullet}^k}$-tight closure is nothing but classical tight closure. That is,  the \textit{tight closure} $I^*$ of an ideal $I \subseteq R$ is defined to be the ideal consisting of all elements $x \in R$ for which there exists $c \in R^{\circ}$ such that $cx^q \in I^{[q]}$ for all large $q=p^e$. 
We say that $R$ is \textit{weakly F-regular} if all ideals in $R$ are tightly closed (i.e.,~$I^*=I$ for all ideals $I \subseteq R$) and that $R$ is \textit{F-regular} if all of its local rings are weakly F-regular. 
The reader is referred to \cite{HH1} for the classical tight closure theory.  
\end{rem}

The test ideal $\tau(R)$ plays a central role in the classical tight closure theory. 
Using the $\a_{\bullet}^k$-tight closure of zero submodule, we define the ideal ${\tau}(\a_{\bullet}^k)$, which is a generalization of  the test ideal $\tau(R)$. 

\begin{propdef}[\textup{cf. \cite[Proposition-Definition 2.9]{Ha}, \cite[Definition-Theorem 6.5]{HY}}]\label{taudef}
Let $\a_{\bullet}=\{\a_m \}$ be a graded family of ideals on an excellent reduced ring $R$ of characteristic $p>0$ and let $k \ge 1$ be an integer. 
Let $E =\bigoplus_{\m} E_R(R/\m)$ be the direct sum, taken over all maximal ideals $\m$ of $R$, of the injective hulls of the residue fields $R/\m$. 

$(1)$ The following ideals are equal to each other. 
\begin{enumerate}
\renewcommand{\labelenumi}{\textup{(\roman{enumi})}}
\item $\displaystyle\bigcap_M \Ann_R(0^{*\a_{\bullet}^{k}}_M)$, where $M$ runs through all finitely generated $R$-modules. 

\item $\displaystyle\bigcap_{M\subset E} \Ann_R(0^{*\a_{\bullet}^{k}}_M)$, where $M$ runs through all finitely generated $R$-submodules of $E$. 
\item $\displaystyle\bigcap_{J\subseteq R} (J:J^{*\a_{\bullet}^{k}})$, where $J$ runs through all ideals of $R$. 
\end{enumerate}
We denote them by ${\tau}(\a_{\bullet}^k)$ $($or $\tau(k \cdot \a_{\bullet}))$ and call it the \textit{generalized test ideal} associated to $\a_{\bullet}$ with exponent $k$. 
Given an ideal $\a \subseteq R$ such that $\a \cap R^{\circ} \ne \emptyset$ and a real number $t>0$, 
if $\a_{\bullet}=\{\a_m\}$ is defined by $\a_m=\a^{\lceil tm \rceil}$, we simply denote this ideal by ${\tau}(\a^t)$. 
Also, just for convenience, we decree that ${\tau}(\a_{\bullet}^{0})={\tau}(R)$.

$(2)$ If $R$ is a $\Q$-Gorenstein normal local ring, then 
$${\tau}(\a_{\bullet}^k)=\Ann_R(0_E^{*\a_{\bullet}^k}).$$

$(3)$ If $(R,\m)$ is a $d$-dimensional Gorenstein local ring and $x_1, \dots, x_d$ is a system of parameters for $R$, then 
$${\tau}(\a_{\bullet}^k)=\bigcap_{l \ge 0}((x_1^l, \dots, x_d^l):(x_1^l, \dots, x_d^l)^{*\a_{\bullet}^k})=\Ann_R(0_{H^d_{\m}(R)}^{*\a_{\bullet}^k}).$$
\end{propdef}

\begin{rem}[\textup{cf.~\cite[Observation 2.8]{Ha}}]
The notation being the same as that in Proposition-Definition \ref{taudef}, 
the generalized test ideal ${\tau}(\a_{\bullet}^k)$ is equal to the unique maximal element among the set of ideals $\{{\tau}(\a_{p^e}^{k/{p^e}})\}_{e \ge 0}$ with respect to inclusion. 
The existence of a maximal element follows from the ascending chain condition on ideals $($because $R$ is Noetherian$)$ and the uniqueness follows from the inclusion ${\tau}(\a_{q}^{k/q}) \subseteq {\tau}(\a_{qq'}^{k/qq'})$ for any powers $q$, $q'$ of $p$. 
If $\a_{\bullet}$ is a descending filtration, then the ideal ${\tau}(\a_{\bullet}^k)$ is equal to the unique maximal element among the set of ideals $\{{\tau}(\a_{m}^{k/{m}})\}_{m \ge 1}$. 
\end{rem}

The notion of $\a_{\bullet}^k$-test elements is useful to study the behavior of the generalized test ideal $\tau(\a_{\bullet}^k)$. 
\begin{defn}[\textup{cf.~\cite[Definition 6.3]{HY}}]\label{testdef}
Let $\a_{\bullet}$ be a graded family of ideals on a ring $R$ of characteristic $p > 0$ and let $k \ge 1$ be an integer. 
An element $d \in R^{\circ}$ is called an {\it $\a_{\bullet}^{k}$-test element} if for every ideal $I \subseteq R$ and every $x \in R$, the following holds: $x \in I^{*\a_{\bullet}^{k}}$ if and only if $dx^q\a_{kq} \subseteq I^{[q]}$ for all powers $q=p^e$ of $p$. 
\end{defn}

An $\a_{\bullet}^k$-test element exists in nearly every ring of interest. 
\begin{prop}[\textup{cf.~\cite[Theorem 1.7]{HY}}]\label{exist}
Let $R$ be a reduced ring of characteristic $p>0$ and let $c \in R^{\circ}$. 
Assume that one of the following conditions holds:
\begin{enumerate}
\item 
$R$ is F-finite and the localized ring $R_c$ is strongly F-regular. 
\item
$R$ is an algebra of finite type over an excellent local ring $B$, 
and the localized ring $R_c$ is Gorenstein and F-regular. 
\end{enumerate}
Then some power $c^n$ of $c$ is an $\a_{\bullet}^k$-test element for all graded family of ideals $\a_{\bullet}$ on $R$ and for all integers $k \ge 1$. 
\end{prop}

\begin{proof}
One can prove this proposition by an argument similar to that in the proof of \cite[Theorem 1.7]{HY}, but we will give a sketch of the proof for the reader's convenience. 
\par 
(1) Take a power $c^n$ which satisfies \cite[Remark 3.2]{HH2}. 
Let $I$ be an ideal of $R$, and fix any $z \in I^{*\a_{\bullet}^k}$ and any power $q$ of $p$. 
Since $z \in I^{*\a_{\bullet}^k}$, there exists $d \in R^{\circ}$ such that 
$dz^Q\a_{kQ} \subseteq I^{[Q]}$ for every power $Q$ of $p$. 
By the choice of $c^n$, there exists a power $q'$ of $p$ and 
$\phi  \in \Hom_R(R^{1/q'},R)$ such that $\phi(d^{1/q'}) = c^n$. 
Since $dz^{qq'}(\a_{kq})^{[q']} \subseteq dz^{qq'}\a_{kqq'} \subseteq I^{[qq']}$, 
one has $d^{1/q'}z^q \a_{kq} R^{1/q'} \subseteq I^{[q]}R^{1/q'}$. 
Applying $\phi$ to both sides gives $c^nz^q \a_{kq} \subseteq I^{[q]}$. 
Hence $c^n$ is an $\a_{\bullet}^k$-test element. 
\par 
(2) Put $R' = R \otimes_B \widehat{B}$, where $\widehat{B}$ denotes the completion of $B$. 
Since $B$ is excellent, $R \to R'$ is faithfully flat with regular fibers.  
In particular, $R'$ is an reduced algebra of finite type over $\widehat{B}$, and 
$R'_c$ is Gorenstein and $F$-regular by \cite[Theorem 7.3 (c)]{HH3}. 
Thus, one can assume that $B$ is complete.   
Then using an $\Gamma$-construction argument (see \cite[Sections 6,7]{HH3}), 
one can reduce the problem to the case (1). 
The reader is referred to \cite[Theorem 6.1, Lemma 6.13, Lemma 6.19]{HH3} for details.
\end{proof}

When the ring is the quotient of a regular local ring, we have a criterion for the triviality of the generalized test ideal $\tau(\a_{\bullet}^k)$. 

\begin{prop}[\textup{cf.~\cite[Theorem 1.12]{Fe}}]\label{Fedder}
Let $(S,\m)$ be a complete regular local ring of characteristic $p>0$ and $I \subsetneq S$ a radical ideal. Let $\a_{\bullet}=\{\a_m \}_{\m \ge 1}$ be a graded family of ideals on $S$ and let $k \ge 1$ be an integer. 
Denote $R=S/I$ and $\a_{R,\bullet}=\{\a_m R\}_{m \ge 1}$.
Fix any $c \in S \setminus I$ whose image in $R$ is an $\a_{R,\bullet}^k$-test element.
Then the generalized test ideal ${\tau}(\a_{R,\bullet}^k)$ is trivial if and only if there exists $q=p^e$ such that $c(I^{[q]}:I)\a_{kq} \not\subseteq \m^{[q]}$.
\end{prop}
\begin{proof}
The proof is similar to those in \cite[Theorem 1.12]{Fe} and \cite[Proposition 2.6]{HW}.  
\end{proof}

In fixed prime characteristic, the generalized test ideal ${\tau}(\a_{\bullet}^k)$ satisfies properties analogous to those of the asymptotic multiplier ideal $\J(\a_{\bullet}^k)$. 

\begin{lem}[\textup{\cite[Lemma 4.5]{Ta}}]\label{inclusiontau}
Let $\a_{\bullet}=\{\a_m \}$ be a graded family of ideals on a ring $R$ of characteristic $p>0$. Then for all integers $k, l \ge 0$, 
$$\a_l {\tau}(\a_{\bullet}^k) \subseteq {\tau}(\a_{\bullet}^{k+l}).$$
\end{lem}
\begin{proof}
It is enough to show that $(0_M^{*\a_{\bullet}^k}:\a_l) \supseteq 0_M^{*\a_{\bullet}^{k+l}}$ for all finitely generated $R$-modules $M$, but it is immediate because $\a_l^{[q]} \a_{kq} \subseteq \a_{(k+l)q}$ for every $q=p^e$. 
\end{proof}

Given graded families of ideals $\a_{\bullet}, \b_{\bullet}$ on a ring $R$ of characteristic $p>0$  
and integers $k, l \ge 1$,  we can define $\a_{\bullet}^k\b_{\bullet}^l$-tight closure as follows:  
if $N$ is a submodule of an $R$-module $M$, then an element $z \in M$ is in the 
$\a_{\bullet}^{k}\b_{\bullet}^{l}$-tight closure $N^{*\a_{\bullet}^{k}\b_{\bullet}^{l}}_M$ of $N$ in $M$ 
if and only if there exists $c \in R^{\circ}$ such that $cz^q\a_{kq}\b_{lq} \subseteq N^{[q]}_M$ for all large $q = p^e$.
The \textit{generalized test ideal} $\tau(\a_{\bullet}^k\b_{\bullet}^l)$ is defined by $\tau(\a_{\bullet}^k\b_{\bullet}^l)=\bigcap_M \Ann_R(0_M^{*\a_{\bullet}^k\b_{\bullet}^l})$, where $M$ runs through all finitely generated $R$-modules. 

\begin{thm}[\textup{\cite[Proposition 2.10]{Ha}, cf.~\cite[Proposition 4.4]{Ta}}]\label{subtau}
Let $R$ be a complete regular local ring of characteristic $p>0$ or an F-finite regular ring of characteristic $p>0$. 
Let $\a_{\bullet}$, $\b_{\bullet}$ be graded families of ideals on $R$ and fix any integers $k, l \ge 0$. Then 
$${\tau}(\a_{\bullet}^{k}\b_{\bullet}^l) \subseteq {\tau}(\a_{\bullet}^{k}){\tau}(\b_{\bullet}^{l}).$$
In particular, 
$${\tau}(\a_{\bullet}^{kl}) \subseteq {\tau}(\a_{\bullet}^{k})^l.$$
\end{thm}

\section{Localization of generalized test ideals}
In \cite{HT}, Hara and the first-named author proved that the formation of generalized test ideals commutes with localization under the assumption of F-finiteness. 
In this section,  we study the behavior of $\a_{\bullet}^k$-tight closure and the ideal $\tau(\a_{\bullet}^k)$ under localization without the assumption of F-finiteness. 

\begin{prop} \label{ideal-loc}
Let $(R,\m)$ be an excellent equidimensional reduced local ring of characteristic $p >0$, 
$\a_{\bullet}$ be a graded family of ideals on $R$ and $k \ge 1$ be an integer. 
Let $W$ denote a multiplicatively closed subset of $R$ 
and $\a_{W,\bullet}$ denote the extension of $\a_{\bullet}$ on $R_W$. 
If $I$ is an ideal of $R$ generated by a subsystem of parameters for $R$, then
\[
 (I R_W)^{*\a_{W,\bullet}^k} = I^{*\a_{\bullet}^k}R_W.
\]
\end{prop}

\begin{proof} 
Since the ring $R$ is excellent and reduced, by Proposition \ref{exist}, it admits an $\b_{\bullet}^l$-test element (say, $d \in R^{\circ}$) for any graded family of ideals $\b_{\bullet}$ on $R$ and for any integer $l \ge 1$.  
We first prove the following claim using an argument similar to that in the proof of \cite[Proposition 2.6]{HH6}. 
\begin{cln}
For any  given $c \in R^{\circ}$, there exists $e_0 \in \N$ 
such that $\xi \in 0^{*\a_{\bullet}^k}_{H_{\m}^d(R)}$ holds 
whenever $\xi \in H_{\m}^d(R)$ and $c\xi^q \a_{kq} =0$ in $H_{\m}^d(R)$ for some $q=p^e \ge p^{e_0}$.
\end{cln}  
\begin{proof}[Proof of Claim 1]
Let $0^F_{H^d_{\m}(R)}$ denote the Frobenius closure of zero in $H^d_{\m}(R)$, that is, $\xi \in H^d_{\m}(R)$ is in $0^F_{H^d_{\m}(R)}$ if and only if $\xi^q=0$ in $H_{\m}^d(R)$ for some power $q=p^e$. 
For each $e \in \N$, we denote 
$$N_e = \{\xi \in H^d_{\m}(R) \mid c\xi^{p^e} \a_{kp^e} \subseteq 0^F_{H^d_{\m}(R)} \}.$$ 
Now let $\xi \in N_{e+1}$ and put $q=p^e$. 
Then $c\xi^{pq}\a_{kpq} \subseteq 0^F_{H^d_{\m}(R)}$ by definition. 
In particular, there exists $q'=p^{e'}$ such that 
$c^{pq'} \xi^{pqq'} \a_{kq}^{[pq']} \subseteq c^{q'} \xi^{pqq'} \a_{kpq}^{[q']}=0$ in $H_{\m}^d(R)$. 
This implies that $c\xi^{q} \a_{kq} \subseteq 0^F_{H^d_{\m}(R)}$, that is, $\xi \in N_e$. 
Hence  $\{N_{e}\}$ is a decreasing sequence on $e$. 
Thus one can choose $e_0 \in \N$ such that $N_e = N_{e_0}$ for all $e \ge e_0$,  because $H_{\m}^d(R)$ is Artinian. 
Then we can easily see that  $N_{e_0}  \subseteq 0_{H^d_{\m}(R)}^{*\a_{\bullet}^k}$, which implies the claim.   
\end{proof}

Let $ x_1,\ldots,x_d$ be a system of parameters for $R$
and we may assume that $I$ is generated by a subsystem of parameters $x_1, \dots , x_h$ $(h \le d)$.
Then we prove the following claim$:$
\begin{cln}
Let $c \in R^{\circ}$ and fix $e_0$ given by Claim 1. 
If $cz^q \a_{kq} \subseteq I^{[q]}$ for some $q=p^e \ge p^{e_0}$, then $z \in I^{*\a_{\bullet}^k}$. 
\end{cln}
\begin{proof}[Proof of Claim 2]
By a standard argument, we can reduce to the case where $h=d$, that is, $I=(x_1, \dots, x_d)$. 
If we put $\xi =[ z + (x_1,\ldots, x_d)] \in H_{\m}^d(R)$, 
then by Claim 1, we have $\xi \in 0^{*\a_{\bullet}^k}_{H_{\m}^d(R)}$,  
because $c\xi^q \a_{kq} =0$ in $H_{\m}^d(R)$. 
Hence $d\xi^Q\a_{kQ} = 0$ in $H^d_{\m}(R)$ for all powers $Q$ of $p$.  
That is, there exists an integer $\ell =\ell(Q) \ge 0$ depending on $Q$ such that 
$d z^Q \a_{kQ} (x_1\cdots x_d)^{\ell} \subseteq (x_1^{\ell+Q},\ldots,x_d^{\ell+Q})$. 
By the colon-capturing property for (classical) tight closure, we get 
$d z^Q \a_{kQ} \subseteq (I^{[Q]})^{*}$ and $d^2 z^Q \a_{kQ} \subseteq I^{[Q]}$ for all powers $Q$. 
This means that $z \in I^{*\a_{\bullet}^k}$, as required. 
\end{proof}

Finally, we prove the assertion of the proposition.   
It is enough to show that $(IR_W)^{*\a_{W,{\bullet}}^k} \subseteq I^{*\a_{\bullet}^k}R_W$.   
Suppose that $\alpha \in (IR_W)^{*\a_{W,{\bullet}}^k}$. 
By definition, there exist $c \in (R_W)^{\circ} \cap R$ and an integer $e_1$ such that 
$c\alpha^q \a_{kq} R_W \subseteq I^{[q]}R_W$ for all $q=p^e \ge p^{e_1}$. 
We may assume that $c \in R^{\circ}$ by prime avoidance. 
By Claim 2, we can take a positive integer $e_2$ such that 
$z \in I^{*\a_{\bullet}^k}$ holds whenever $cz^q\a_{kq} \subseteq I^{[q]}$ 
for some $q=p^e \ge p^{e_2}$. 
Fix any $e \ge e_3:=\max\{e_1,e_2\}$ and put $q=p^e$. 
Choose $u \in W$ such that $uc\alpha^q \a_{kq} \subseteq I^{[q]}$. 
Then $c(u \alpha)^q \a_{kq} \subseteq I^{[q]}$ and thus  $u\alpha \in I^{*\a_{\bullet}^k}$. 
That is, $\alpha \in I^{*\a_{\bullet}^k}R_W$, as required. 
\end{proof}

\begin{cor}[\textup{cf. \cite[Proposition 3.1]{HT}}]\label{tau-loc}
Let $(R,\m)$ be a complete Gorenstein reduced local ring of characteristic $p >0$, $\a_{\bullet}$ be a graded family of ideals on $R$ and $k \ge 1$ be an integer. 
Let $W$ denote a multiplicatively closed subset of $R$ and $\a_{W,\bullet}$ denote the extension of $\a_{\bullet}$ on $R_W$. Then  
\[
{\tau}(\a_{W, \bullet}^k) = {\tau}(\a_{\bullet}^k)R_W.
\]
\end{cor}

\begin{proof}
The proof is based on arguments similar to those in \cite{Sm2}.
First, we will show that $\tau(\a_{\bullet}^k)R_W \subseteq \tau(\a_{W,\bullet}^k)$. 
To see this, let $c \in \tau(\a_{\bullet}^k)$ and $z \in (IR_W)^{*\a_{W,\bullet}^k}\cap R$, 
where $I$ is any given ideal of $R$ such that $IR_W \ne R_W$. 
Let $P$ is a prime ideal of $R$ such that $P \cap W = \emptyset$ and $\a_{P,\bullet}$ denote the extension of $\a_{\bullet}$ on $R_P$. 
Take an $R$-sequence $x_1,\ldots,x_h$ in $P$ whose images form a system of parameters for $R_{P}$. 
By Proposition \ref{ideal-loc}, 
$$\tau(\a_{\bullet}^k)((x_1^l, \dots, x_h^l)R_{P})^{*\a_{P,\bullet}^k}
= \tau(\a_{\bullet}^k)(x_1^l, \dots, x_h^l)^{*\a_{\bullet}^k}R_{P}
\subseteq (x_1^l, \dots, x_h^l)R_{P}$$ for every $l \in \N$. 
Thus $\tau(\a_{\bullet}^k)R_{P} \subseteq  \tau(\a_{P,\bullet}^k)$, because $R$ is Gorenstein. 
Since the image $\frac{z}{1}$ lies in $(IR_{P})^{*\a_{P,\bullet}^k}$, we have that $\frac{~cz~}{1} = \frac{~c~}{1} \cdot \frac{~z~}{1} \in IR_{P}$ for all such primes $P$. 
It then follows that $\frac{cz}{1} \in IR_W$, that is, $\frac{c}{1} \in \tau(\a_{W,\bullet}^k)$. 

Next, we will prove the converse.
We start with the case where $W = R \setminus P$ and $P \subset R$ is a prime ideal of height $h$. 
Then 
\[
 \tau(\a_{P,\bullet}^k)=\Ann_{R_{P}} 0^{*\a_{P,\bullet}^k}_{H^h_{P R_{P}}(R_{P})}
\]
by definition. 
On the other hand, since $R$ is complete and $H_{\m}^d(R)=E$ is Artinian, 
\[
 \tau(\a_{\bullet}^k) R_{P} = (\Ann_R N)R_{P} = \Ann_{R_{P}} N^{\vee_{\m}\vee_{P}},
\]
where $N=0^{*\a_{\bullet}^k}_{H_{\m}^d(R)}$ and 
$$N^{\vee_{\m}\vee_{P}} = \Hom_R(\Hom_R(N,H^d_{\m}(R)),H^h_{P R_{P}}(R_{P})) \subseteq H^h_{P R_{P}}(R_{P}),$$ by \cite[Lemma 3.1(iii)]{Sm1}. 
So it is enough to show the following claim in order to prove 
$\tau(\a_{\bullet}^k)R_{P} \supseteq \tau(\a_{P,\bullet}^k):$
\begin{cl}
$N^{\vee_{\m}\vee_{P}} \subseteq 0^{*\a_{P,\bullet}^k}_{H^h_{P R_{P}}(R_{P})}$ in $H^h_{P R_{P}}(R_{P})$. 
\end{cl}
\begin{proof}[Proof of Claim]
Take  a system of parameters $x_1,\ldots,x_d$ for $R$ such that $\frac{x_1}{1},\ldots.\frac{x_h}{1}$ forms a system of parameters for $R_{P}$. 
Suppose that 
$\eta = \left[ \frac{~z~}{1} + \left(\frac{x_1^t}{1},\cdots, \frac{x_h^t}{1} \right) \right] \in H^h_{P R_{P}}(R_{P})$ belongs to $N^{\vee_{\m}\vee_{P}}$.
To see that $\a_{kq} \eta^q \subseteq N^{\vee_{\m}\vee_{P}}$ for all $q=p^e$,  
we may assume that $t=1$ without loss of generality. 
Since $N^{\vee_{\m}\vee_{P}} = \Ann_{H^h_{P R_{P}}(R_{P})} \tau(\a_{\bullet}^k)R_{P}$ by \cite[Lemma 2.1 (iv)]{Sm2}, 
we have  
$\frac{~z~}{1} \tau(\a_{\bullet}^k)R_{P} \subseteq \left(\frac{x_1}{1},\cdots, \frac{x_h}{1} \right)R_{P}$. 
Take any $a \in R \setminus P$ such that $az\tau(\a_{\bullet}^k) \subseteq (x_1,\ldots,x_h)$ and put 
$\eta_{\ell} = \left[az + (x_1,\ldots,x_h,x_{h+1}^\ell,\ldots,x_d^\ell) \right] \in H_{\m}^d(R)$ for each integer $l \ge 1$. 
Then since $\eta_\ell \in \Ann_{H_{\m}^d(R)} \tau(\a_{\bullet}^k) = 0^{*\a_{\bullet}^k}_{H_{\m}^d(R)}$ (by Matlis dual),  
we have that $\a_{kq} \eta_\ell^q \subseteq \Ann_{H_{\m}^d(R)} \tau(\a_{\bullet}^k)$ for all $q=p^e$ by an argument similar to that in the proof of \cite[Proposition 1.15]{HY}.
That is, for every $\ell \ge 1$, 
$(az)^q \a_{kq} \tau(\a_{\bullet}^k) \subseteq (x_1^q,\ldots,x_h^q, x_{h+1}^{\ell q},\ldots,x_d^{\ell q})$. 
In particular, 
$\left(\frac{z}{1}\right)^q \a_{kq} \tau(\a_{\bullet}^k)R_P \subseteq \left(\frac{x_1^q}{1},\ldots, \frac{x_h^q}{1}\right)R_{P}$ 
and thus $\a_{kq} \eta^q \subseteq N^{\vee_{\m}\vee_{P}}$ for all $q=p^e$. 

Note that there exists some $\frac{c}{1} \in \tau(\a_{\bullet}^k)R_{P} \cap (R_{P})^{\circ}$. 
Then $\frac{c}{1}$ kills every element of $N^{\vee_{\m}\vee_{P}}$.
This means that for any $\eta \in N^{\vee_{\m}\vee_{P}}$, one has $\frac{c}{1}\a_{kq} \eta^q=0$ in $H^h_{P R_{P}}(R_{P})$ for all $q=p^e$. 
We conclude that $N^{\vee_{\m}\vee_{P}} \subseteq 0^{*\a_{P,\bullet}^k}_{H^h_{P R_{P}}(R_{P})}$. 
\end{proof}

Finally, we consider the general case. 
Suppose that $c \in \tau(\a_{W,\bullet}^k) \cap R$. 
To see  $c \in \tau(\a_{\bullet}^k)R_W \cap R$, it suffices to show that 
$c \in \tau(\a_{\bullet}^k)R_{P} \cap R$ for every prime ideal $P \subset R$ such that $P \cap W = \emptyset$.  
Take a system of parameters $x_1,\ldots,x_d$ such that $\frac{x_1}{1},\ldots,\frac{x_h}{1}$ forms a 
system of parameters for $R_{P}$. 
Then by the definition of $c$ and Proposition \ref{ideal-loc}, 
$$c((x_1^{\ell}, \dots, x_h^{\ell})R_{P})^{*\a_{P,\bullet}^k} = c ((x_1^{\ell}, \dots, x_h^{\ell})R_W)^{*\a_{W,\bullet}^k}R_{P} \subseteq (x_1^{\ell}, \dots, x_h^{\ell}) R_{P}$$
for every $\ell \ge 1$. 
Thus, $c \in \tau(\a_{P,\bullet}^k) \cap R = \tau(\a_{\bullet}^k)R_{P} \cap R$. 
\end{proof}

Making use of the above corollary, we can show a Skoda-type theorem for symbolic powers of ideals.  
\begin{prop}[\textup{cf. \cite[Theorem 2.12]{Ha}}]\label{skodatau}
Let $R$ be an excellent Gorenstein reduced local ring of characteristic $p>0$ and let $I$ be an ideal of $R$ such that $I \cap R^{\circ} \ne \emptyset$. 
Suppose that the residue field of each of the rings $R_P$ is infinite when $P$ is an associated prime of $I$.  
Let $I^{(\bullet)}=\{I^{(m)}\}$ denote the graded family of symbolic powers of $I$ and $h$ denote the largest analytic spread of $IR_P$ as $P$ runs through the associated primes of $I$. 
Then for every integer $k \ge 0$,  
$${\tau}((h+k) \cdot I^{(\bullet)}) \subseteq I^{(k+1)}.$$
\end{prop}

\begin{proof}
Let $I_P^{(\bullet)}$ denote the extension of $I^{(\bullet)}$ on $R_P$. 
Since after localization at $P$ the symbolic and ordinary powers of $I$ are the same,  $I_P^{(\bullet)}=\{(IR_P)^m\}_{m \ge1}$. 
By assumption, $I R_P$ has a reduction ideal generated by at most $h$ elements for each associated prime $P$ of $I$.  
Then, by \cite[Theorem 2.1]{HY} and the first half of the proof of Corollary \ref{tau-loc}, 
\begin{align*}
{\tau}((h+k) \cdot I^{(\bullet)})R_P \subseteq {\tau}((h+k) \cdot I_P^{(\bullet)})&={\tau}((IR_P)^{h+k})\\
&\subseteq (IR_P)^{k+1}
\end{align*}
for every associated prime $P$ of $I$. 
Thus one has ${\tau}((h+k) \cdot I^{(\bullet)}) \subseteq I^{(k+1)}$, as required. 
\end{proof}

\section{Symbolic powers in positive characteristic}
Eisenbud-Mazur asked in \cite{EM} a question on the behavior of symbolic square of ideals in regular local rings of equal characteristic zero.  
When the ring is of positive characteristic or of mixed characteristic, counterexamples to their question are known (see \cite{EM} for the positive characteristic case and \cite{KR} for the mixed characteristic case). 
Nevertheless, Hochster-Huneke proved in \cite{HH5} analogous results to their question in positive characteristic, using the classical tight closure theory. 
In this section, using generalized test ideals, we give a slight generalization of Hochster-Huneke's results. 

\begin{thm}[\textup{cf. \cite[Theorem 3.5]{HH5}}]\label{charpEM}
Let $(R,\m)$ be an excellent regular local ring of characteristic $p>0$ and $I \subsetneq R$ be any nonzero ideal.
Let $h$ be the largest analytic spread of $IR_P$ as $P$ runs through the associated primes of $I$. 
Then for all integers $n \ge 1$ and $k \ge 0$, one has 
$$I^{(hn+kn+1)}\subseteq \m (I^{(k+1)})^n.$$
In particular if $h=2$, then 
$$I^{(3)} \subseteq \m I.$$ 
\end{thm}

\begin{proof}
We may assume that $R$ is a complete local ring: suppose that we have the assertion in the complete regular case. 
Although $\widehat{I}=I\widehat{R}$ may have more associated primes, the biggest analytic spread as one localizes at these cannot increase (see \cite[Discussion 2.3 (c)]{HH4} for details). Also, $I^{(m)}\widehat{R}=\widehat{I}^{(m)}$ for every integer $m \ge 1$. Thus,
$$I^{(hn+kn+1)}= \widehat{I}^{(hn+kn+1)} \cap R \subseteq \widehat{\m} (\widehat{I}^{(k+1)})^n \cap R=\m (I^{(k+1)})^n,$$
as required, by the faithful flatness of the completion. 

Moreover, if $R/\m$ is a finite field (in this case, $R$ is F-finite), then we replace $R$ by $R[t]_M$ and $I$ by $IR[t]_M$ where $t$ is an indeterminate and $M$ is a maximal ideal of $R[t]$ generated by $\m$ and $t$. 
Note that the issues are unaffected by this replacement, because the associated primes of $IR[t]_M$ is simply those of the form $PR[t]_M$ where $P$ is an associated prime of $I$ (see \cite[Discussion 2.3 (b)]{HH4} for further explanation). 
By this trick, we can assume that $R$ is a complete regular local ring with infinite residue field or an F-finite regular local ring such that the residue field of each of the rings $R_P$ is infinite when $P$ is an associated prime of $I$.  

Let $I^{(\bullet)}=\{I^{(m)}\}$ denote the graded family of symbolic powers of $I$ and $l$ denote the largest integer such that ${\tau}(l \cdot I^{(\bullet)})=R$. 
Such an integer $l$ always exists, because ${\tau}(0 \cdot P^{(\bullet)})=R$. 
Note that the ideal $\tau((l+1) \cdot P^{(\bullet)})$ is contained in $\m$. 
By Lemma \ref{inclusiontau}, Theorem \ref{subtau} and Proposition \ref{skodatau}, one has 
\begin{align*}
I^{(hn+kn+1)}=I^{(hn+kn+1)}\tau(l \cdot I^{(\bullet)}) & \subseteq \tau((hn+kn+l+1) \cdot I^{(\bullet)})\\
&\subseteq \tau((l+1) \cdot I^{(\bullet)}) \tau((h+k) \cdot I^{(\bullet)})^n\\
& \subseteq \m (I^{(k+1)})^n.
\end{align*}
\end{proof}

\begin{lem}\label{inversioncharp}
Let $R$ be an excellent regular local ring of characteristic $p>0$ and let $I \subsetneq R$ be an ideal of height at least two such that $R/I$ is F-pure. 
Let $I^{(\bullet)}=\{I^{(m)}\}_{m \ge 0}$ denote the graded family of symbolic powers of $I$ and $\widehat{I}^{(\bullet)}=\{\widehat{I}^{(m)}\}_{m \ge 0}$ denote the extension of  $I^{(\bullet)}$ on the completion $\widehat{R}$. 
Then the generalized test ideal ${\tau}(\widehat{R},  \widehat{I}^{(\bullet)})$ associated to $\widehat{I}^{(\bullet)}$ is trivial. 
Moreover, if $R$ is F-finite, then the generalized test ideal ${\tau}(R, I^{(\bullet)})$ associated to $I^{(\bullet)}$ is also  trivial. 
\end{lem}

\begin{proof}
By \cite[Theorem 1.12]{Fe},  the ring $R/I$ is F-pure if and only if $(I^{[q]}:I) \not\subset \m^{[q]}$ for all $q=p^e$.
Also, by Proposition \ref{Fedder},  the generalized test ideal $\tau(\widehat{R},  \widehat{I}^{(\bullet)})$ is trivial if and only if $I^{(q)}\widehat{R} =\widehat{I}^{(q)} \not\subset (\m\widehat{R})^{[q]}$ for some $q=p^e$. 
Therefore, it is enough to show that $(I^{[q]}:I) \subseteq I^{(q)}$ for every $q=p^e$. 
Let $P_1, \dots, P_k$ be the minimal prime ideals of the radical ideal $I$. Then, by definition, $I^{(q)}=\bigcap_{i=1}^kP_i^{(q)}$. On the other hand, we have the following claim. 

\begin{cl}
$(I^{[q]}:I)= \bigcap_{i=1}^k(P_i^{[q]}:P_i)$
\end{cl}
\begin{proof}[Proof of Claim]
It is enough to show that $(P_i^{[q]}:P_i)=(P_i^{[q]}:I)$ for all $i=1, \dots, k$. Since the inclusion $(P_i^{[q]}:P_i)\subseteq (P_i^{[q]}:I)$ is clear, we will prove the reverse inclusion.
Let $x \in (P_i^{[q]}:I)$. Then $xP_1 \cdots P_k \subseteq P_i^{[q]}$ and it follows that $xP_iR_{P_i} \subseteq P_i^{[q]}R_{P_i}$. Since the Frobenius map is flat, $P_i^{[q]}$ is a $P_i$-primary ideal. This implies that $xP_i \in P_i^{[q]}$. 
\end{proof}
Thus, we may assume that $I$ is a prime ideal of height $h\ge 2$. 
Since the Frobenius map is flat, $I^{[q]}$ is $I$-primary and then $(I^{[q]}:I)$ is also $I$-primary. 
Therefore, one has 
\begin{align*}
(I^{[q]}:I)&=(IR_I^{[q]}:IR_I) \cap R\\
&=(I^{h(q-1)}R_I+I^{[q]}R_I) \cap R\\
&\subseteq I^{(q)},
\end{align*}
because $R_I$ is a regular local ring of dimension $h \ge 2$. 

The latter assertion is immediate, because the formation of generalized test ideals commutes with completion if the ring is F-finite (this follows from an argument similar to that of the proof of \cite[Proposition 3.2]{HT}). 
\end{proof}

\begin{thm}[\textup{cf. \cite[Theorem 3.6]{HH5}}]\label{Fpurecase}
Let $R$ be an excellent regular ring of characteristic $p>0$ and let $I \subsetneq R$ be an ideal of height at least two such that $R/I$ is F-pure.  
Let $h$ denote the largest height of any minimal prime of $I$. 
Then for all integers $n \ge 1$ and $k \ge 0$, one has 
$$I^{(hn+kn-1)} \subseteq (I^{(k+1)})^{n}.$$ 
In particular if $h=2$, then 
$$I^{(3)} \subseteq I^2.$$ 
\end{thm}

\begin{proof}
The problem reduces to the local case, and then by essentially the same argument as that in the proof of Theorem \ref{charpEM}, we can assume that $R$ is a complete regular local ring with infinite residue field or an F-finite regular local ring such that the residue field of each of the rings $R_P$ is infinite when $P$ is an associated prime of $I$. 
By virtue of Lemma \ref{inversioncharp}, the generalized test ideal $\tau(1 \cdot I^{(\bullet)})$ is trivial. 
Then for all integers $n \ge 1$ and $k \ge 0$, it follows from Lemma \ref{inclusiontau},Theorem \ref{subtau} and Proposition \ref{skodatau} that 
\begin{align*}
I^{(hn+kn-1)}=I^{(hn+kn-1)}\tau(1 \cdot I^{(\bullet)}) &\subseteq \tau((hn+kn) \cdot I^{(\bullet)})\\
&\subseteq \tau((h+k) \cdot I^{(\bullet)})^n\\
&\subseteq (I^{(k+1)})^n.
\end{align*} 
\end{proof}

\begin{rem}
The proof of Theorem \ref{Fpurecase} tells us the following: given an ideal $I$ of an excellent regular ring $R$ of characteristic $p>0$ such that $I \cap R^{\circ} \ne \emptyset$, if the generalized test ideal $\tau(l \cdot I^{(\bullet)})$ associated to $I^{(\bullet)}$ is trivial for some integer $l \ge 1$, then one has $I^{(hn+kn-l)} \subseteq (I^{(k+1)})^{n}$ for all integers $n \ge 1$ and $k \ge 0$, where $h$ is the largest analytic spread of $IR_P$ as $P$ runs through the associated primes of $I$.  
\end{rem}

In general, Eisenbud-Mazur's conjecture fails in positive characteristic (see \cite{EM}), but there is no counterexample to the codimension two case. 
Theorem \ref{Fpurecase} suggests that symbolic powers of an ideal in a regular ring behave nicely in case the quotient of the ring by the ideal is F-pure. 
We expect that even in positive characteristic, their conjecture holds true for ideals of codimension two if the quotient of the regular ring by the ideal is F-pure. 

\begin{conj}
Let $(R,\m)$ be an excellent regular local ring of characteristic $p>0$ and let $I \subsetneq R$ be an unmixed ideal of height two such that $R/I$ is F-pure.  Then $I^{(2)} \subseteq \m I$? 
\end{conj}

This conjecture is known to hold when $R/I$ is Cohen-Macaulay (because such $I$ is a licci ideal and then $P^{(2)} \subseteq \m P$ is established in \cite{EM}). 
Also, it is easy to check that the conjecture holds for squarefree monomial ideals. 

\section{Symbolic powers in characteristic zero}

Using the standard descent theory of \cite[Chapter 2]{HH7}, we generalize the positive characteristic results in previous section to the equal characteristic case. 
\begin{thm}
Let $(R,\m)$ be a regular local ring essentially of finite type over a field of characteristic zero and $I \subsetneq R$ be any nonzero ideal.
Let $h$ be the largest analytic spread of $IR_P$ as $P$ runs through the associated primes of $I$. 
Then for all integers $n \ge 1$ and $k \ge 0$, one has 
$$I^{(hn+kn+1)}\subseteq \m(I^{(k+1)})^n.$$
In particular if $h=2$, then 
$$I^{(3)} \subseteq \m I.$$
\end{thm}

\begin{proof}
We employ the same strategy as that of the proof of Theorem \ref{charpEM} and use asymptotic multiplier ideals instead of generalized test ideals.  
Let $I^{(\bullet)}=\{I^{(m)}\}$ denote the graded family of symbolic powers of $I$ and $l$ denote the largest integer such that ${\J}(l \cdot I^{(\bullet)})=R$. 
Note that the ideal $\J((l+1) \cdot I^{(\bullet)})$ is contained in $\m$. 
Applying \cite[Theorem 11.1.19]{La},  the subadditivity formula (\cite[Theorem 11.2.3]{La}) and Skoda's theorem (\cite[Theorem 9.6.21]{La}), one has 
\begin{align*}
I^{(hn+kn+1)}=I^{(hn+kn+1)}\J(l \cdot I^{(\bullet)}) & \subseteq \J((hn+kn+l+1) \cdot I^{(\bullet)})\\
&\subseteq \J((l+1) \cdot I^{(\bullet)}) \J((h+k) \cdot I^{(\bullet)})^n\\
& \subseteq \m (I^{(k+1)})^n.
\end{align*}
\end{proof}

\begin{defn}
Let $R$ be a ring which is finitely generated over a field $k$ of characteristic zero. 
Then the ring $R$ is said to be of \textit{dense F-pure type} if there exist a finitely generated $\Z$-algebra $A \subseteq k$ and a finitely generated $A$-algebra $R_A$ which is free over $A$ such that $R \cong R_A \otimes_A k$ and that for all maximal ideals $\mu$ in a Zariski dense subset of $\Spec A$ with residue field $\kappa=A/\mu$, the fiber rings $R_A \otimes_A \kappa$ are F-pure. 
\end{defn}

\begin{thm}
Let $R$ be a regular algebra of essentially of finite type over a field of characteristic zero and let $I \subsetneq R$ be an ideal of height at least two such that $R/I$ is of dense F-pure type.  
Let $h$ denote the largest height of any minimal prime of $I$. 
 Then for all integers $n \ge 1$ and $k \ge 0$, one has
 $$I^{(hn+kn-1)} \subseteq (I^{(k+1)})^n.$$ 
 In particular if $h=2$, then 
 $$I^{(3)} \subseteq I^2.$$ 
\end{thm}

\begin{proof} 
We use the standard descent theory of \cite[Chapter 2]{HH7} to reduce the problem to the positive characteristic case. 
This reduction step is essentially same as the one used in the proof of \cite[Theorem 4.4]{HH4} and \cite[Theorem 4.2]{HH5}. 
After reduction to characteristic $p>0$, the assertion immediately follows from Theorem \ref{Fpurecase}. 
\end{proof}

\end{document}